\documentclass[12pt,oneside,reqno]{amsart}
\usepackage{graphicx}
\usepackage{mathrsfs}
\usepackage{stmaryrd}
\usepackage{amsfonts}
\usepackage{enumerate,amsmath,amssymb,amsthm}
\usepackage{color}
\newcommand{\dif}{\mathrm{d}}

\newcommand{\be}{\begin{eqnarray}}
\newcommand{\ee}{\end{eqnarray}}
\newcommand{\ce}{\begin{eqnarray*}}
\newcommand{\de}{\end{eqnarray*}}
\newtheorem{theorem}{Theorem}[section]
\newtheorem{lemma}[theorem]{Lemma}
\newtheorem{remark}[theorem]{Remark}
\newtheorem{definition}[theorem]{Definition}
\newtheorem{proposition}[theorem]{Proposition}
\newtheorem{Examples}[theorem]{Examples}
\newtheorem{corollary}[theorem]{Corollary}

\def\[{{\Big[}}
\def\]{{\Big]}}
\def\<{{\langle}}
\def\>{{\rangle}}
\def\({{\Big(}}
\def\){{\Big)}}

\def\no{\nonumber}
\def\bt{\begin{theorem}}
\def\et{\end{theorem}}
\def\bl{\begin{lemma}}
\def\el{\end{lemma}}
\def\br{\begin{remark}}
\def\er{\end{remark}}
\def\bx{\begin{Examples}}
\def\ex{\end{Examples}}
\def\bd{\begin{definition}}
\def\ed{\end{definition}}
\def\bp{\begin{proposition}}
\def\ep{\end{proposition}}
\def\bc{\begin{corollary}}
\def\ec{\end{corollary}}

\def\cK{{\mathcal K}}

\def\cM{{\mathcal M}}

\def\cP{{\mathcal P}}

\def\mE{{\mathbb E}}

\def\mN{{\mathbb N}}

\def\mP{{\mathbb P}}

\def\mR{{\mathbb R}}

\def\sB{{\mathscr B}}
\def\sC{{\mathscr C}}

\def\sF{{\mathscr F}}

\def\geq{\geqslant}
\def\leq{\leqslant}

\begin{document}

\allowdisplaybreaks

\title{Stability for Stochastic McKean-Vlasov Equations with Non-Lipschitz Coefficients}

\author{Xiaojie Ding and Huijie Qiao$^*$}

\thanks{{\it AMS Subject Classification(2000):} 60H15}

\thanks{{\it Keywords:} Stochastic McKean-Vlasov equations, exponential stability of moments, exponentially 2-ultimate boundedness, almost surely asymptotic stability.}

\thanks{This work was supported by NSF of China (No. 11001051, 11371352) and China Scholarship Council under Grant No. 201906095034.}

\thanks{$*$ Corresponding author: Huijie Qiao}

\subjclass{}

\date{}
\dedicatory{School of Mathematics,
Southeast University\\
Nanjing, Jiangsu 211189, P.R.China\\
Email: 220171469@seu.edu.cn, hjqiaogean@seu.edu.cn}

\begin{abstract}
In this paper we consider the stability for a type of stochastic McKean-Vlasov equations with non-Lipschitz coefficients. First, sufficient conditions are given for the exponential stability of the second moments for their solutions in terms of a Lyapunov function. Then we weaken the conditions and furthermore obtain exponentially 2-ultimate boundedness of their solutions. After this, the almost surely asymptotic stability of their solutions is proved. Finally we give an example to motivate the choice of Lyapunov functions.
\end{abstract}

\maketitle \rm

\section{Introduction}

Given a complete filtered probability space $(\Omega,\mathscr{F},\{\mathscr{F}_t\}_{t\in[0,\infty)},\mP)$. Consider the following stochastic McKean-Vlasov equation on $\mR^d$:
\be\left\{\begin{array}{ll}
X_t=\xi+\int_0^tb(X_s,\mu_s)ds+\int_0^t\sigma(X_s,\mu_s)d W_s,\\
\mu_s=\mP_{X_s}= $the probability  distribution of~$X_s,
\end{array}
\label{eq1}
\right.
\ee
where $\xi$ is a $\sF_0$-measurable random variable with $\mE|\xi|^{4}<\infty$, $W_{\cdot}=(W_{\cdot}^1,W_{\cdot}^2,\cdots,W_{\cdot}^l)$ is a $(\sF_t)_{t\geq0}$-adapted standard $l$-dimensional Brownian motion and the coefficients $b:\mR^d\times\cM_{\lambda^2}(\mR^d)\rightarrow{\mR^d}, \sigma:\mR^d\times\cM_{\lambda^2}(\mR^d)\rightarrow{\mR^d}\times{\mR^l}$ are Borel measurable. ($\cM_{\lambda^2}(\mR^d)$ is defined in Section \ref{nn})

If $b, \sigma$ are independent of $\mu_s$,  Eq.(\ref{eq1}) becomes a usual stochastic differential equation(SDE). Moreover, in recent years, the stability for solutions of SDEs has been studied extensively in the literature. Most of these papers are concerned with exponential stability of $p$-th moments for their solutions, exponential stability of sample paths for their solutions and exponential stability, exponentially 2-ultimate boundedness, or almost surely asymptotic stability of their solutions. Let us mention some works. For linear SDEs, Arnold collected a number of results on exponential stability of their solutions in his monograph \cite{la}. For SDEs in infinite dimensional Hilbert spaces, Ichikawa \cite{IA} proved the stability of moments and exponential stability of sample paths for their solutions under Lipschitz and linear growth conditions. For SDEs with jumps, Deng, Krsti\'c and Williams \cite{DW} studied the almost surely asymptotic stability by using a strong Markov property. Later, the second named author and Duan \cite{QD} offered some general conditions of exponentially 2-ultimate boundedness of solutions for SDEs with jumps. For SDEs with jumps in infinite dimensional Hilbert spaces, Bao, Truman and Yuan \cite{BTY} discussed the almost surely asymptotic stability for their solutions under local Lipschitz condition but without a linear growth condition.

If $b, \sigma$ depend on $\mu_s$,  Eq.(\ref{eq1}) is called as a stochastic McKean-Vlasov equation(SMVE). And there are few results about the stability of its solution due to its specialty including distributions. Recently, for a semilinear stochastic McKean-Vlasov evolution equation, Govindan and Ahmed \cite{GA} investigated the exponential stability for its solution under Lipschitz and linear growth conditions.

In this paper, we study the stability of SMVEs under non-Lipschitz conditions. In \cite{DH}, we have proved that Eq.$(\ref{eq1})$ has a unique strong solution. Here we continue and consider three types of stability for the strong solution to Eq.$(\ref{eq1})$. First, we give sufficient conditions to prove the exponential stability of the second moment in terms of Lyapunov functions. Then by a similar way, it is shown that the exponentially 2-ultimate boundedness of its solution holds. Finally, motivated by \cite{BTY}, we take a nonrandom initial condition and prove the almost surely asymptotic stability under more general conditions of Lyapunov functions.

The rest of the paper is organized as follows. In Section \ref{fram}, we recall some basic notation, give some necessary assumptions and extend the classical It\^o's formula to the It\^o formula for SMVEs. And then we prove the exponential stability of the second moment for the strong solution to Eq.$(\ref{eq1})$ in Section \ref{esom}. In Section \ref{eub}, the exponentially 2-ultimate boundedness of the strong solution for Eq.$(\ref{eq1})$ is investigated. Then, the almost surely asymptotic stability of the strong solution for Eq.$(\ref{eq1})$ is proved in Section \ref{asas}. In section \ref{eg}, an example is given to explain our results.

The following convention will be used throughout the paper: $C$ with or without indices will denote
different positive constants whose values may change from one place to another.

\section{The Framework}\label{fram}

In the section, we recall some basic notation, give some necessary concepts and assumptions and extend the classical It\^o's formula to the It\^o formula for SMVEs.

\subsection{Notation}\label{nn}

In the subsection, we introduce notation used in the sequel.

Let $\sB(\mR^d)$ be the Borel $\sigma$-algebra on $\mR^d$ and $\cP({\mR^d})$ be the space of all probability measures defined on $\sB(\mR^d)$ carrying the usual topology of weak convergence. Let $C(\mR^d)$ be the collection of continuous functions on $\mR^d$. For convenience, we shall use $\mid\cdot\mid$ and $\parallel\cdot\parallel$  for norms of vectors and matrices, respectively. Furthermore, let $\langle\cdot$ , $\cdot\rangle$ denote the scalar product in $\mR^d$. Let $A^*$ denote the transpose of the matrix $A$.

Define the Banach space
$$
C_\rho(\mR^d):=\left\{{\varphi\in{C(\mR^d)},\|{\varphi}\|_{C_\rho(\mR^d)}
=\sup_{x\in{\mR^d}}\frac{|{\varphi(x)}|}{(1+|{x}|)^2}+\sup_{x\neq{y}}\frac{|{\varphi(x)-\varphi(y)}|}{|{x-y}|}<\infty}\right\}.
$$
Let $\cM_{\lambda^2}^s(\mR^d)$ be the Banach space of signed measures $m$
on $\sB(\mR^d)$ satisfying
\ce
\|m\|_{\lambda^2}^2:=\int_{\mR^d}(1+|{x}|)^2\,|m|(\dif x)<\infty,
\de
where $|m|=m^{+}+m^{-}$ and $m=m^{+}-m^{-}$ is the Jordan decomposition of $m$. Let
$\cM_{\lambda^2}(\mR^d)=\cM_{\lambda^2}^s(\mR^d)\bigcap\cP(\mR^d)$ be the set of probability
measures on $\sB(\mR^d)$ with finite second order moments. We put on $\cM_{\lambda^2}(\mR^d)$ a topology induced by the
following metric:
\ce
\rho(\mu,\nu):=\sup_{\parallel{\varphi}\parallel_{C_\rho(\mR^d)\leq1}}\left|{\int_{\mR^d}\varphi(x)\mu(dx)-\int_{\mR^d}\varphi(x)\nu(dx)}\right|.
 \de
Then $(\cM_{\lambda^2}(\mR^d),\rho)$ is a complete metric space.

Next, we recall the definition of the derivative for $f:\cM_{\lambda^2}(\mR^d)\rightarrow\mR$ with respect to a probability measure as introduced in his lecture by Lions \cite{Lion}. A function $f:\cM_{\lambda^2}(\mR^d)\rightarrow\mR$ is differential at $\mu\in \cM_{\lambda^2}(\mR^d)$, if for $\tilde{f}(\xi):=f(\mP_\xi),\xi\in L^2(\Omega,\mathscr{F},\mP;\mR^d)$, there exists some $\zeta\in L^2(\Omega,\mathscr{F},\mP;\mR^d)$ with $\mP_\zeta=\mu$ such that $\tilde{f}$ is Fr\'echet differentiable at $\zeta$, that is, there exists a linear continuous mapping $D\tilde{f}(\zeta):L^2(\Omega,\mathscr{F},\mP;\mR^d)\rightarrow\mR$ such that for any $\eta\in L^2(\Omega,\mathscr{F},\mP;\mR^d)$
\ce
\tilde{f}(\zeta+\eta)-\tilde{f}(\zeta)=D\tilde{f}(\zeta)(\eta)+o(|\eta|_{L^2}),  \quad|\eta|_{L^2}\rightarrow0.
\de
Since $D\tilde{f}(\zeta)\in L(L^2(\Omega,\mathscr{F},\mP;\mR^d),\mR)$, from the Riesz representation theorem, there exists a $\mP$-a.s. unique variable $\vartheta\in L^2(\Omega,\mathscr{F},\mP;\mR^d)$ such that for all $\eta\in L^2(\Omega,\mathscr{F},\mP;\mR^d)$
\ce
D\tilde{f}(\zeta)(\eta)=(\vartheta,\eta)_{L^2}=\mE[\vartheta\cdot\eta].
\de
Lions \cite{Lion} proved that there is a Borel measurable function $h:\mR^d\rightarrow\mR^d$ which depends on the distribution $\mP_\zeta$ rather than $\zeta$ itself, such that $\vartheta=h(\zeta)$. Therefore, for $\xi\in L^2(\Omega,\mathscr{F},\mP;\mR^d)$
\ce
f(\mP_\xi)-f(\mP_\zeta)=\mE[h(\zeta)(\xi-\zeta)]+o(|\xi-\zeta|_{L^2}).
\de
We call $\partial_\mu f(\mP_\zeta)(y):=h(y),y\in \mR^d$, the derivative of $f:\cM_{\lambda^2}(\mR^d)\rightarrow\mR$ at $\mP_\zeta,\zeta\in L^2(\Omega,\mathscr{F},\mP;\mR^d)$.

\bd\label{f1}
We say that $f\in C^1(\cM_{\lambda^2}(\mR^d))$, if there exists for all $\xi\in L^2(\Omega,\mathscr{F},\mP;\mR^d)$ a $\mP_{\xi}$-modification of $\partial_\mu f(\mP_\xi)(\cdot)$, again denoted by $\partial_\mu f(\mP_\xi)(\cdot)$, such that $\partial_\mu f:\cM_{\lambda^2}(\mR^d)\times\mR^d\rightarrow\mR^d$ is continuous, and we identify this continuous function $\partial_\mu f$ as the derivative of $f$.
\ed

\bd\label{r1}
The function $f$ is said to be in $C_b^{1;1}(\cM_{\lambda^2}(\mR^d))$, if $f\in C^1(\cM_{\lambda^2}(\mR^d))$ and $\partial_\mu f$ is bounded and Lipschitz continuous, that is, there exists a real constant $C>0$ such that
\ce
&(i)&|\partial_\mu f(\mu)(x)|\leq C,\mu\in\cM_{\lambda^2}(\mR^d),x\in\mR^d;\\
&(ii)&|\partial_\mu f(\mu)(x)-\partial_\mu f(\mu')(x')|\leq C(\rho(\mu,\mu')+|x-x'|),\mu,\mu'\in\cM_{\lambda^2}(\mR^d),x,x'\in\mR^d.
\de
\ed

\bd\label{f2}
We say that $f\in C^2(\cM_{\lambda^2}(\mR^d))$, if for any $\mu\in \cM_{\lambda^2}(\mR^d)$, $f\in C^1(\cM_{\lambda^2}(\mR^d))$ and $\partial_\mu f(\mP_\xi)(\cdot)$ is differentiable, and its derivative $\partial_y\partial_\mu f:\cM_{\lambda^2}(\mR^d)\times\mR^d\rightarrow\mR^d\otimes\mR^d$ is continuous.
\ed

\bd\label{r2}
The function $f$ is said to be in $C_b^{2;1}(\cM_{\lambda^2}(\mR^d))$, if $f\in C^2(\cM_{\lambda^2}(\mR^d))\cap C_b^{1;1}(\cM_{\lambda^2}(\mR^d))$ and its derivative $\partial_y\partial_\mu f$ is bounded and Lipschitz continuous.
\ed

\bd\label{r2}
The function $\Phi$ is said to be in $C_b^{2,2;1}(\mR^d\times\cM_{\lambda^2}(\mR^d))$, if \\
(i) $\Phi$ is bi-continuous in $(x,\mu)$;\\
(ii) For any $x$, $\Phi(x,\cdot)\in C_b^{2;1}(\cM_{\lambda^2}(\mR^d))$, and for any $\mu\in\cM_{\lambda^2}(\mR^d)$, $\Phi(\cdot, \mu)\in C^2(\mR^d)$.\\
If $\Phi\in C_b^{2,2;1}(\mR^d\times\cM_{\lambda^2}(\mR^d))$ and $\Phi\geq 0$, we say that $\Phi\in C_{b,+}^{2,2;1}(\mR^d\times\cM_{\lambda^2}(\mR^d))$.
\ed

\bd
The function $\Phi$ is said to be in $\sC(\mR^d\times\cM_{\lambda^2}(\mR^d))$, if $\Phi\in C^{2,2}(\mR^d\times\cM_{\lambda^2}(\mR^d))$ and for any compact set $\cK\subset\mR^d\times\cM_{\lambda^2}(\mR^d)$,
$$
\sup\limits_{(x,\mu)\in\cK}\int_{\mR^d}\left(\|\partial_y\partial_\mu \Phi(x,\mu)(y)\|^2+|\partial_\mu \Phi(x,\mu)(y)|^2\right)\mu(\dif y)<\infty.
$$
If $\Phi\in\sC(\mR^d\times\cM_{\lambda^2}(\mR^d))$ and $\Phi\geq 0$, we say that $\Phi\in\sC_+(\mR^d\times\cM_{\lambda^2}(\mR^d))$.
\ed

\subsection{Some assumptions}\label{ass}

In the subsection, we give out some assumptions.
\begin{enumerate}[($\bf{H}_{1.1}$)]
\item The functions $b, \sigma$ are continuous in $(x,\mu)$ and satisfy for $(x,\mu)\in\mR^{d}\times{\cM_{\lambda^2}(\mR^d)}$
\ce
{|{b(x,\mu)}|}^2+\|\sigma(x,\mu)\|^2\leq{L_1(1+|{x}|^2+\|{\mu}\|^2_{\lambda^2})},
\de
where $L_1>0$ is a constant.
\end{enumerate}

\begin{enumerate}[($\bf{H}_{1.2}$)]
\item The functions $b,\sigma$ satisfy for $(x_1,\mu_1), (x_2,\mu_2)\in\mR^{d}\times{\cM_{\lambda^2}(\mR^d)}$
\ce
&&2\langle{x_1-x_2,b(x_1,\mu_1)-b(x_2,\mu_2)}\rangle+\parallel{\sigma(x_1,\mu_1)-\sigma(x_2,\mu_2)}\parallel^2\\
&\leq&{L_2\(\kappa_1(|x_1-x_2|^2)+\kappa_2\left(\rho^2(\mu_1,\mu_2)\right)\)},
\de
where $L_2>0$ is a constant, and $\kappa_i(x), i=1, 2$ are two positive, strictly increasing, continuous concave function and satisfy $\kappa_i(0)=0$,  $\int_{0^+}\frac{1}{\kappa_1(x)+\kappa_2(x)}dx=\infty$.
\end{enumerate}

By \cite[Theorem 3.1]{DH}, we know that Eq.$(\ref{eq1})$ has a unique strong solution denoted as $X_t$ under $(\bf{H}_{1.1})$-$(\bf{H}_{1.2})$. And then we assume some other conditions to prove the exponential stability of the second moment for $X_t$.

\begin{enumerate}[($\bf{H}_{2.1}$)]
\item There exists a function $v:\mR^d\times\cM_{\lambda^2}(\mR^d)\rightarrow\mR$ satisfying\\
$(i)$ $v\in \sC_{+}(\mR^d\times\cM_{\lambda^2}(\mR^d))$,\\
$(ii)$
$$
\int_{\mR^d}\(\mathscr{L^\mu}v(x,\mu)+\alpha v(x,\mu)\)\mu(dx)\leq 0, 
$$
where $\mathscr{L^\mu}$ is defined as
\ce
\(\mathscr{L^\mu}v\)(x,\mu)&:=&\(b^i\partial_{x_i} v\)(x,\mu)+\frac{1}{2}\((\sigma\sigma^*)^{ij}\partial_{x_ix_j}^2 v\)(x,\mu)
+\int_{\mR^d}b^i(y,\mu)(\partial_\mu v)_i(x,\mu)(y)\mu(dy)\\
&&+\int_{\mR^d}\frac{1}{2}(\sigma\sigma^*)^{ij}(y,\mu)\partial_{y_i}(\partial_\mu v)_j(x,\mu)(y)\mu(dy),
\de
and $\alpha>0$ is a constant, \\
$(iii)$ 
$$
a_1\int_{\mR^d}|x|^2\mu(dx)\leq\int_{\mR} v(x,\mu)\mu(dx)\leq a_2\int_{\mR^d}|x|^2\mu(dx),
$$ 
where $a_1, a_2>0$ are two constants.
\end{enumerate}

Here and hereafter we use the convention that the repeated indices stand for the summation. In the following, we weaken ($\bf{H}_{2.1}$) to show the exponentially 2-ultimate boundedness. 

\begin{enumerate}[($\bf{H}_{2.2}$)]
\item There exists a function $v:\mR^d\times\cM_{\lambda^2}(\mR^d)\rightarrow\mR$ satisfying\\
$(i)$ $v\in\sC(\mR^d\times\cM_{\lambda^2}(\mR^d))$,\\
$(ii)$
$$
\int_{\mR^d}\(\mathscr{L^\mu}v(x,\mu)+\alpha v(x,\mu)\)\mu(dx)\leq M_1,
$$
$(iii)$ 
$$
a_1\int_{\mR^d}|x|^2\mu(dx)-M_2\leq \int_{\mR} v(x,\mu)\mu(dx)\leq a_2\int_{\mR^d}|x|^2\mu(dx)+M_3,
$$
where $M_1$, $M_2$, $M_3\geq0$ are constants.
\end{enumerate}

Next, we strengthen ($\bf{H}_{2.1}$) to prove the almost surely asymptotic stability of $X_t$. Here we introduce a function class. Let $\varSigma$ denote the family of functions $\gamma:\mR_+\rightarrow\mR_+$, which are continuous, strictly increasing, and $\gamma(0)=0$. And $\varSigma_\infty$ means the family of functions $\gamma\in \varSigma$ with $\gamma(x)\rightarrow\infty$ as $x\rightarrow\infty$.

\begin{enumerate}[($\bf{H}^{\prime}_{1.1}$)]
\item The function $b$ is continuous in $(x,\mu)$ and satisfies for $(x,\mu)\in\mR^{d}\times{\cM_{\lambda^2}(\mR^d)}$
\ce
{|{b(x,\mu)}|}^2\leq{L'_1(1+|{x}|^2+\|{\mu}\|^2_{\lambda^2})},
\de
where $L'_1>0$ is a constant, and $\sigma$ is bounded.
\end{enumerate}

\begin{enumerate}[($\bf{H}_{2.3}$)]
\item There exists a function $v:\mR^d\times\cM_{\lambda^2}(\mR^d)\rightarrow\mR$ satisfying\\
$(i)$ $v\in C_{b,+}^{2,2;1}(\mR^d\times\cM_{\lambda^2}(\mR^d))$,\\
$(ii)$
$\mathscr{L^\mu}v(x,\mu)+\alpha v(x,\mu)\leq 0$,\\
$(iii)$
$\gamma_1(|x|)\leq v(x,\mu)\leq \gamma_2(|x|)$,
where $\gamma_i\in\varSigma_\infty(i=1,2).$
\end{enumerate}

\br
(i) ($\bf{H}^{\prime}_{1.1}$) is stronger than ($\bf{H}_{1.1}$) and ($\bf{H}_{2.3}$) is stronger than ($\bf{H}_{2.1}$).

(ii) Note that when $b, \sigma, v$ do not depend on the distribution, the operator $\mathscr{L^\mu}$ reduces to
\ce
\mathscr{L^\mu}v(x)=\mathscr{L}v(x):={\(b^i\partial_{x_i}}v\)(x)+\frac{1}{2}\((\sigma\sigma^*)^{ij}\partial_{x_ix_j}^2v\)(x).
\de
And then ($\bf{H}_{2.1}$) ($\bf{H}_{2.2}$) become classical conditions and have appeared in \cite{la}. 
\er

\subsection{It\^{o}'s formula for SMVEs}\label{SP}
Next we will extend the classical It\^{o}'s formula to SDEs depending on the distribution $\mu$.
\bp\label{if}
Let $f:\mR^d\times\cM_{\lambda^2}(\mR^d)\rightarrow\mR$ be such that $f\in \sC(\mR^d\times\cM_{\lambda^2}(\mR^d))$. Then, under $(\bf{H}_{1.1})$-$(\bf{H}_{1.2})$, the following It\^{o}'s formula holds:
\be
f(X_t,\mu_t)-f(X_s,\mu_s)&=&\int_s^t{(b^i\partial_{x_i}}f)(X_u,\mu_u)du+\frac{1}{2}\int_s^t\((\sigma\sigma^*)^{ij}\partial_{x_ix_j}^2f\)(X_u,\mu_u)du\no\\
&&+\int_s^t\int_{\mR^d}b^i(y,\mu_u)(\partial_\mu f)_i(X_u,\mu_u)(y)\mu_u(\dif y)du\no\\
&&+\frac{1}{2}\int_s^t\int_{\mR^d}(\sigma\sigma^*)^{ij}(y,\mu_u)\partial_{y_i}(\partial_\mu f)_j(X_u,\mu_u)(y)\mu_u(\dif y)du\no\\
&&+\int_s^t(\sigma^{ij}\partial_i f)(X_u,\mu_u)dW_u^j, \quad 0\leq s<t.
\label{ito1}
\ee
\ep
\begin{proof}
Under $(\bf{H}_{1.1})$-$(\bf{H}_{1.2})$, we know that $X_{\cdot}$ satisfies Eq.(\ref{eq1}), i.e.
\ce
X_t=\xi+\int_0^tb(X_s,\mu_s)ds+\int_0^t\sigma(X_s,\mu_s)d W_s.
\de
And then by the H\"older inequality and the isometry formula it holds that for any $T>0$ and $0\leq t\leq T$, 
\be
\mE|X_t|^2&\leq& 3\mE|\xi|^2+3\mE\left|\int_0^tb(X_s,\mu_s)ds\right|^2+3\mE\left|\int_0^t\sigma(X_s,\mu_s)d W_s\right|^2\no\\
&\leq&3\mE|\xi|^2+3T\mE\int_0^t\left|b(X_s,\mu_s)\right|^2ds+3\mE\int_0^t\left\|\sigma(X_s,\mu_s)\right\|^2ds\no\\
&\leq&3\mE|\xi|^2+3(T+1)\mE\int_0^tL_1(1+|X_s|^2+\|\mu_s\|^2_{\lambda^2})ds\no\\
&\leq&3\mE|\xi|^2+3(T+1)\int_0^tL_1(1+\mE|X_s|^2+2\mE(1+|X_s|^2))ds\no\\
&\leq&3\mE|\xi|^2+9(T+1)TL_1+9(T+1)L_1\int_0^t\mE|X_s|^2ds,
\label{2dedu}
\ee
where we use $\mE|\xi|^4<\infty$ and the fact that $\|\mu_s\|^2_{\lambda^2}\leq\mE(1+|X_s|)^2\leq 2\mE(1+|X_s|^2)$. The Gronwall inequality admits us to obtain that 
\be
\sup\limits_{t\in[0,T]}\mE|X_t|^2\leq C.
\label{2semo}
\ee
By the similar deduction to the above (\ref{2dedu}), one can get that 
\ce
\mE|X_t|^4&\leq& 27\mE|\xi|^4+81L_1^2(T^{7/4}+T^{3/2})+81L_1^2(T^{3/4}+T^{1/2})\int_0^t\|\mu_s\|^4_{\lambda^2}ds\\
&&+81L_1^2(T^{3/4}+T^{1/2})\int_0^t\mE|X_s|^4ds\\
&\overset{(\ref{2semo})}\leq&C+81L_1^2(T^{3/4}+T^{1/2})\int_0^t\mE|X_s|^4ds.
\de
From the Gronwall inequality, it follows that
\be
\sup\limits_{t\in[0,T]}\mE|X_t|^4\leq C.
\label{4semo}
\ee
So, by $(\bf{H}_{1.1})$ we know that 
\ce
\mE\int_0^T(|b(X_r,\mu_r)|^2+\|\sigma(X_r,\mu_r)\|^4)\dif r&\leq& L_1\mE\int_0^T\left(1+|X_r|^2+\|\mu_r\|^2_{\lambda^2}\right)\dif r\\
&&+3L^2_1\mE\int_0^T\left(1+|X_r|^4+\|\mu_r\|^4_{\lambda^2}\right)\dif r\\
&\overset{(\ref{2semo})(\ref{4semo})}<&\infty,
\de
Then by \cite[Proposition A.8]{WDL}, we know that (\ref{ito1}) holds. The proof is complete.
\end{proof}

\section{The exponential stability of the second moment}\label{esom}

In the section, we study the exponential stability of the second moment for the strong solution to Eq.$(\ref{eq1})$.

\bt\label{sm}
Assume that $(\bf{H}_{1.1})$-$(\bf{H}_{1.2})$ and $(\bf{H}_{2.1})$ hold. Then $X_t$ satisfies
\ce
\mE|X_t|^2\leq \frac{a_2}{a_1} e^{-\alpha t}\mE|\xi|^2, \quad\quad t\geq0.
\de
\et
\begin{proof}
Applying It\^{o}'s formula to $e^{\alpha t}v(X_t,\mu_t)$, we have
\ce
&&e^{\alpha t}v(X_t,\mu_t)-v(\xi,\mu_0)\no\\
&=&\int_0^t\alpha e^{\alpha s}v(X_s,\mu_s)ds+\int_0^te^{\alpha s}(b^i\partial_{x_i}v)(X_s,\mu_s)ds+\frac{1}{2}\int_0^t e^{\alpha s}\((\sigma^*\sigma)^{ij}\partial_{x_ix_j}^2v\)(X_s,\mu_s)ds\no\\
&&+\int_0^t\int_{\mR^d}e^{\alpha s}b^i(y,\mu_s)(\partial_\mu v)_i(X_s,\mu_s)(y)\mu_s(dy)ds+\int_0^te^{\alpha s}(\partial_iv\sigma^{ij})(X_s,\mu_s)dW_s^j\no\\
&&+\frac{1}{2}\int_0^t\int_{\mR^d}e^{\alpha s}(\sigma^*\sigma)^{ij}(y,\mu_s)\partial_{y_i}(\partial_\mu v)_j(X_s,\mu_s)(y)\mu_s(dy)ds\no\\
&=&\int_0^te^{\alpha s}\[\alpha v(X_s,\mu_s)+(b^i\partial_{x_i}v)(X_s,\mu_s)+\frac{1}{2}\((\sigma^*\sigma)^{ij}\partial_{x_ix_j}^2v\)(X_s,\mu_s)\no\\
&&+\int_{\mR^d}b^i(y,\mu_s)(\partial_\mu v)_i(X_s,\mu_s)(y)\mu_s(dy)+\frac{1}{2}\int_{\mR^d}(\sigma^*\sigma)^{ij}(y,\mu_s)\partial_{y_i}(\partial_\mu v)_j(X_s,\mu_s)(y)\mu_s(dy)\]ds\no\\
&&+\int_0^te^{\alpha s}(\partial_iv\sigma^{ij})(X_s,\mu_s)dW_s^j\no\\
&=&\int_0^te^{\alpha s}\(\alpha v(X_s,\mu_s)+\mathscr{L^\mu}v(X_s,\mu_s)\)ds+\int_0^te^{\alpha s}(\partial_iv\sigma^{ij})(X_s,\mu_s)dW_s^j.
\de
Localizing and taking the expectation on both sides of the above equality, by the Fatou Lemma one can get that
\ce
e^{\alpha t}\mE v(X_t,\mu_t)-\mE v(\xi,\mP_\xi)\leq\mE\[\int_0^te^{\alpha s}\(\alpha v(X_s,\mu_s)+\mathscr{L^\mu}v(X_s,\mu_s)\)ds\].
\de
Then it follows from $(ii)$ in $(\bf{H}_{2.1})$ that
\ce
e^{\alpha t}\mE v(X_t,\mu_t)-\mE v(\xi,\mP_\xi)
\leq\int_0^t e^{\alpha s}\mE\(\alpha v(X_s,\mu_s)+\mathscr{L^\mu}v(X_s,\mu_s)\)ds\leq 0,
\de
which yields
\ce
\mE v(X_t,\mu_t)\leq e^{-\alpha t}\mE v(\xi,\mP_\xi).
\de
Moreover, by $(iii)$ in $(\bf{H}_{2.1})$ it holds that
\ce
a_1\mE|X_t|^2\leq\mE v(X_t,\mu_t)\leq e^{-\alpha t}\mE v(\xi,\mP_\xi)\leq a_2e^{-\alpha t}\mE|\xi|^2.
\de
Thus, we obtain
\ce
\mE|X_t|^2\leq \frac{a_2}{a_1} e^{-\alpha t}\mE|\xi|^2.
\de
The proof is complete.
\end{proof}

\section{The exponentially 2-ultimate boundedness}\label{eub}

In the section, we study the exponentially 2-ultimate boundedness for the solution of Eq.$(\ref{eq1})$. First of all, we introduce the exponentially 2-ultimate boundedness.

\bd\label{d2}
If there exist positive constants $K$, $\beta$, $M$ such that
\ce
\mE|X_t|^2\leq Ke^{-\beta t}\mE|\xi|^2+M, \quad t\geq0,
\de
then the solution $X_t$ for Eq.$(\ref{eq1})$ is called exponentially 2-ultimately bounded.
\ed

\bt\label{sp}
Suppose that $(\bf{H}_{1.1})$-$(\bf{H}_{1.2})$ and $(\bf{H}_{2.2})$ hold. Then $X_t$ is exponentially 2-ultimately bounded, i.e.
\ce
\mE|X_t|^2\leq\frac{a_2}{a_1}e^{-\alpha t}\mE|\xi|^2+\frac{\alpha(M_2+M_3)+M_1}{\alpha a_1}, \quad t\geq0.
\de
\et

Since the proof of the above theorem is similar to that in Theorem \ref{sm}, we omit it.

\section{The almost surely asymptotic stability}\label{asas}

In the section, we require that $\xi=x_0$ is non-random and study the almost surely asymptotic stability of the strong solution for Eq.$(\ref{eq1})$.

First of all, we introduce the almost surely asymptotic stability.

\bd\label{d1}
The solution of Eq.$(\ref{eq1})$ is said to be almost surely asymptotically stable if for all $x_0\in\mR^d$, it holds that
\ce
\mP\left\{\lim_{t\rightarrow\infty}|X_t|=0\right\}=1.
\de
\ed

\bt\label{as}
Assume that $(\bf{H}^{\prime}_{1.1})$ $(\bf{H}_{1.2})$ and $(\bf{H}_{2.3})$ hold. Then $X_t$ is almost surely asymptotically stable.
\et
\begin{proof}
To prove that for all $x_0\in\mR^d$,
\ce
\mP\left\{\lim_{t\rightarrow\infty}|X_t|=0\right\}=1,
\de
by $(iii)$ of $(\bf{H}_{2.3})$ we only need to show that for all $x_0\in\mR^d$,
\be
\mP\left\{\lim_{t\rightarrow\infty}v(X_t,\mu_t)=0\right\}=1.
\label{4.2}
\ee

{\bf Step 1.} Assume that $x_0=0$. We prove (\ref{4.2}).

Applying It\^{o}'s formula to $v(X_t,\mu_t)$, one can obtain
\be
v(X_t,\mu_t)&=&v(X_s,\mu_s)+\int_s^t\((b^i\partial_{x_i}v)(X_u,\mu_u)+\frac{1}{2}((\sigma\sigma^*)^{ij}\partial_{x_ix_j}^2v)(X_u,\mu_u)\)du\no\\
&&+\int_s^t\int_{\mR^d}b^i(y,\mu_u)(\partial_\mu v)_i(X_u,\mu_u)(y)\mu_u(dy)du\no\\
&&+\frac{1}{2}\int_s^t\int_{\mR^d}(\sigma\sigma^*)^{ij}(y,\mu_u)\partial_{y_i}(\partial_\mu v)_j(X_u,\mu_u)(y)\mu_u(dy)du\no\\
&&+\int_s^t(\sigma^{ij}\partial_i v)(X_u,\mu_u)dW_u^j,\quad\quad\quad 0\leq s<t.
\label{itofor}
\ee
Thus by $(ii)$ in $(\bf{H}_{2.3})$, it holds that
\ce
v(X_t,\mu_t)&\leq&v(X_s,\mu_s)-\alpha\int_s^tv(X_u,\mu_u)du+\int_s^t(\sigma^{ij}\partial_iv)(X_u,\mu_u)dW_u^j\\
&\leq&v(X_s,\mu_s)+\int_s^t(\sigma^{ij}\partial_iv)(X_u,\mu_u)dW_u^j,  \quad\quad\quad 0\leq s<t.
\de
Set $\mathscr{G}_t:=\sigma(\mathscr{F}_t^W\cup\mathscr{N})$ for $t\geq 0$, where $\{\mathscr{F}_t^W\}_{t\geq 0}$ is the $\sigma$-algebra generated by $W_{\cdot}$ and $\mathscr{N}$ is the collection of all the $\mP$ null sets. And then we have that $v(X_\cdot,\mu_\cdot)$ is adapted to $\{\mathscr{G}_t\}_{t\geq 0}$ when $X_\cdot$ is the strong solution of Eq.$(\ref{eq1})$, and $\int_s^t(\sigma^{ij}\partial_iv)(X_u,\mu_u)dW_u^j$ is a martingale with respect to $\{\mathscr{G}_t\}_{t\geq s}$. Therefore, we obtain
\ce
\mE(v(X_t,\mu_t)\mid\mathscr{G}_s)&\leq&\mE(v(X_s,\mu_s)\mid\mathscr{G}_s)+\mE\left(\int_s^t(\sigma^{ij}\partial_iv)(X_u,\mu_u)dW_u^j\mid\mathscr{G}_s\right)\\
&\leq&v(X_s,\mu_s),
\de
which implies that $v(X_\cdot,\mu_\cdot)$ is a supermartingale with respect to $\{\mathscr{G}_t\}_{t\geq 0}$. If $s=0$ and $X_0=x_0=0$, by $(iii)$ of $(\bf{H}_{2.3})$ and the supermartingale property of $v(X_\cdot,\mu_\cdot)$, it holds that $v(X_t,\mu_t)=0$, a.s. for $t\geq0$. Thus, (\ref{4.2}) is right.

{\bf Step 2.} Assume $x_0\neq0$. We prove (\ref{4.2}). 

Set
\ce
A_1&:=&\{\omega:\liminf_{t\rightarrow\infty} v(X_t,\mu_t)>0\},\\
A_2&:=&\{\omega:\limsup_{t\rightarrow\infty} v(X_t,\mu_t)>0\},
\de
and then we just need to prove $\mP(A_1)=\mP(A_2)=0.$  Let $s=0$, taking the expectation on two sides of (\ref{itofor}), one can have that
\ce
\mE(v(X_t,\mu_t))=v(x_0,\delta_{x_0})+\mE\left(\int_0^t\mathscr{L^\mu}v(X_s,\mu_s)ds\right)=v(x_0,\delta_{x_0})+\int_0^t\mE\mathscr{L^\mu}v(X_s,\mu_s)ds,
\de
where $\delta_{x_0}$ is the Dirac measure in $x_0$. Thus, $(ii)$ and $(iii)$ in $(\bf{H}_{2.3})$ admits us to obtain that
\ce
0\leq\mE\gamma_1(|X_t|)\leq\mE(v(X_t,\mu_t))\leq v(x_0,\delta_{x_0})-\alpha\int_0^t\mE v(X_s,\mu_s)ds, \quad t\geq0,
\de
and furthermore
\ce
\int_0^t\mE v(X_s,\mu_s)ds\leq\frac{v(x_0,\delta_{x_0})}{\alpha}, \quad t\geq0.
\de
Let $t\rightarrow\infty$, we obtain that
\be
\int_0^\infty \mE v(X_s,\mu_s)ds\leq\frac{v(x_0,\delta_{x_0})}{\alpha}<\infty.\label{4.3}
\ee
Hence, by the Fatou lemma it holds that
\ce
\liminf_{t\rightarrow\infty} v(X_t,\mu_t)=0,\quad a.s.
\de
that is, $\mP(A_1)=0$.

Next, we prove that $\mP(A_2)=0$. It follows from the simple calculation that
$$
\mP(A_2)=\mP(A_2A^c_1)+\mP(A_2A_1)=\mP(A_2A^c_1),
$$
where $A^c_1$ stands for the complementary set of $A_1$. Therefore, we only need to prove $\mP(A_2A^c_1)=0$. Note that 
$$
A_2A^c_1=\left\{\omega: \liminf\limits_{t\rightarrow\infty}v(X_t,\mu_t)=0, \limsup\limits_{t\rightarrow\infty} v(X_t,\mu_t)>0\right\}. 
$$
Suppose that $\mP(A_2A^c_1)\neq0$. And then there exists $\varepsilon_1, \varepsilon_2>0$ such that
\be
\mP\left\{v(X_.,\mu_.)\ crosses\ from\ below\ \varepsilon_1\ to\ above\ 2\varepsilon_1\ and\
 back\ infinitely\ many\ times\right\}\geq\varepsilon_2.\label{4.4}
\ee

Here, to make the following deduction convenient, we rewrite Eq.(\ref{eq1}) as
\be
X_t=x_0+\int_0^t\tilde{b}(u, X_u)\dif u+\int_0^t\tilde{\sigma}(u, X_u)\dif W_u, \quad t\geq 0,
\label{Eq2}
\ee
where $\tilde{b}(u, X_u):=b(X_u,\mu_u), \tilde{\sigma}(u, X_u):=\sigma(X_u,\mu_u)$. Since Eq.(\ref{eq1}) has a unique strong solution $(X_t)_{t\geq 0}$ with the initial distribution $\delta_{x_0}$, the distribution family $\{\mu_t\}_{t\geq 0}$ of $(X_t)_{t\geq 0}$ is known. That is, Eq.(\ref{Eq2}) is a nonhomogeneous classical SDE.

Set $\tau_n:=\inf\{t\geq0,|X_t|>n\}.$
Now applying It\^{o}'s formula to $|X_{s\wedge\tau_n}-x_0|^2$ for $s\geq0$, one can get that
\ce
|X_{s\wedge\tau_n}-x_0|^2&=&2\sum\limits_{i=1}^d\int_0^{s\wedge\tau_n}{( X^i_u-x^i_0) \tilde{b}^i(u, X_u)}du+\sum\limits_{i=1}^d\int_0^{s\wedge\tau_n}(\tilde{\sigma}\tilde{\sigma}^*(u, X_u))^{ii}du\\
&&+2\sum\limits_{i,j=1}^d\int_0^{s\wedge\tau_n}( X^i_u-x^i_0)\tilde{\sigma}^{ij}(u, X_u)dW_u^j\\
&\leq&2\int_0^{s\wedge\tau_n}{|X_u-x_0| | \tilde{b}(u, X_u)|}du+\int_0^{s\wedge\tau_n}\|\tilde{\sigma}(u, X_u)\|^2du\\
&&+2\sum\limits_{i,j=1}^d\int_0^{s\wedge\tau_n}( X^i_u-x^i_0) \tilde{\sigma}^{ij}(u, X_u)dW_u^j\\
&\leq&\int_0^{s\wedge\tau_n}|X_u-x_0|^2du+\int_0^{s\wedge\tau_n}|\tilde{b}(u, X_u)|^2du+\int_0^{s\wedge\tau_n}\|\tilde{\sigma}(u, X_u)\|^2du\\
&&+2\sum\limits_{i,j=1}^d\int_0^{s\wedge\tau_n}( X^i_u-x^i_0)\tilde{\sigma}^{ij}(u, X_u)dW_u^j.\\
\de
By $(\bf{H}^{\prime}_{1.1})$ and the Burkholder-Davis-Gundy inequality, we can derive that
\ce
\mE\left(\sup_{0\leq s\leq t}|X_{s\wedge\tau_n}-x_0|^2\right)&\leq&\mE\int_0^{t\wedge\tau_n}|X_u-x_0|^2du+C\mE\int_0^{t\wedge\tau_n}(1+| X_u|^2+n^2)du\\
&&+C\mE\left(\int_0^{t\wedge\tau_n}|X_u-x_0|^2\|\tilde{\sigma}(u, X_u)\|^2du\right)^{1/2}\\
&\leq&\mE\int_0^{t\wedge\tau_n}|X_u-x_0|^2du+C\mE\int_0^{t\wedge\tau_n}(1+|X_u|^2+n^2)du\\
&&+\frac{1}{4}\mE\left(\sup_{0\leq s\leq t}|X_{s\wedge\tau_n}-x_0|^2\right)+C\mE\left(\int_0^{t\wedge\tau_n}\|\tilde{\sigma}(u, X_u)\|^2du\right),
\de
which yields
\ce
\mE\left(\sup_{0\leq s\leq t}|X_{s\wedge\tau_n}-x_0|^2\right)&\leq& C\mE\int_0^{t\wedge\tau_n}|X_u-x_0|^2du+C\mE\int_0^{t\wedge\tau_n}(1+| X_u|^2+n^2)du\\
&&+C\mE({t\wedge\tau_n}).
\de
It follows from the boundedness of $X_t$ that
\ce
\mE\left(\sup_{0\leq s\leq t}|X_{s\wedge\tau_n}-x_0|^2\right)\leq C\mE({t\wedge\tau_n})\leq Ct,
\de
where $C>0$ is depending on $L_1$, $x_0$ and $n$. Then for any $\lambda>0$, Chebyshev's inequality gives that
\be
\mP\left\{\sup_{0\leq s\leq t}|X_{s\wedge\tau_n}-x_0|>\lambda\right\}\leq\frac{Ct}{\lambda^2}.\label{4.5}
\ee

Besides, on one hand, for any $0<\epsilon_1\leq \frac{1}{2}$, we can choose $\Gamma(\cdot)\in\varSigma_\infty$ such that $\frac{v(x_0,\delta_{x_0})}{\Gamma(v(x_0,\delta_{x_0}))}\leq\epsilon_1$. On the other hand, by $(iii)$ of $(\bf{H}_{2.3})$, we know that $\sup\limits_{0\leq s\leq t}\mid X_s\mid\geq\eta(\mid x_0\mid)$ implies $\sup\limits_{0\leq s\leq t}v(X_s,\mu_s)\geq\Gamma(v(x_0,\delta_{x_0}))$, where $\eta(\cdot):=\gamma_1^{-1}\circ\Gamma\circ\gamma_2(\cdot)$ and $\gamma_1^{-1}$ is the inverse function of $\gamma_1$. Thus, we obtain
\be
\mP\left\{\sup_{0\leq s\leq t}\mid X_s\mid\geq\eta(\mid x_0\mid)\right\}&\leq&\mP\left\{\sup\limits_{0\leq s\leq t}v(X_s,\mu_s)\geq\Gamma(v(x_0,\delta_{x_0}))\right\}\no\\
&\leq&\frac{\mE\left(\sup\limits_{0\leq s\leq t}v(X_s,\mu_s)\right)}{\Gamma(v(x_0,\delta_{x_0}))}\no\\
&\leq&\frac{v(x_0,\delta_{x_0})}{\Gamma(v(x_0,\delta_{x_0}))}\leq\epsilon_1, \quad\quad \forall t\geq0.\label{4.1}
\ee

Next, since $v(x,\mu)$ is bicontinuous in $(x,\mu)$, it must be uniformly continuous in $B:=\{(x,\mu)\in\mR^d\times\cM_{\lambda^2}(\mR^d): |x|<\eta(|x_0|), \|\mu\|_{\lambda^2}<\sqrt{2+2\eta^2(|x_0|)}\}$.
Therefore, for any $\epsilon_2>0$, we can choose a function $\gamma\in\varSigma$ such that if for any $(x_1, \mu_1), (x_2, \mu_2)\in B$, $| x-y|<\gamma(\epsilon_2)$ and $\rho(\mu_1,\mu_2)<\gamma(\epsilon_2)$, then
\be
|v(x_1,\mu_1)-v(x_2,\mu_2)|\leq\epsilon_2.
\label{es0}
\ee
Thus, combining (\ref{4.5}) $(\ref{4.1})$ with (\ref{es0}), one can obtain
\ce
&&\mP\left\{\sup_{0\leq s\leq t}|v(X_s,\mu_s)-v(x_0,\delta_{x_0})|>\epsilon_2\right\}\\
&=&\mP\left\{\sup_{0\leq s\leq t}|v(X_s,\mu_s)-v(x_0,\delta_{x_0})|>\epsilon_2,\ \sup_{0\leq s\leq t}|X_s|<\eta(|x_0|)\right\}\\
&&+\mP\left\{\sup_{0\leq s\leq t}|v(X_s,\mu_s)-v(x_0,\delta_{x_0})|>\epsilon_2,\ \sup_{0\leq s\leq t}|X_s|\geq\eta(|x_0|)\right\}\\
&\overset{(\ref{es0})}{\leq}&\mP\left\{\sup_{0\leq s\leq t}|X_s-x_0|>\gamma(\epsilon_2),\ \sup_{0\leq s\leq t}|X_s|<\eta(|x_0|)\right\}\\
&&+\mP\left\{\sup_{0\leq s\leq t}|X_s|\geq\eta(|x_0|)\right\}\\
&\leq&\mP\left\{\sup_{0\leq s\leq t}|X_{s\wedge\tau_{\eta(|x_0|)}}-x_0|>\gamma(\epsilon_2)\right\}+\mP\left\{\sup_{0\leq s\leq t}|X_s|\geq\eta(| x_0|)\right\}\\
&\overset{(\ref{4.5})(\ref{4.1})}{\leq}&\frac{Ct}{\gamma^2(\epsilon_2)}+\epsilon_1.
\de
Taking $t^*$ such that $\frac{Ct^*}{\gamma^2(\epsilon_2)}\leq\frac{1}{4}$ and noting $\epsilon_1\leq\frac{1}{2}$, we get that
\ce
\mP\left\{\sup_{0\leq s\leq t^*}|v(X_s,\mu_s)-v(x_0,\delta_{x_0})|>\epsilon_2\right\}\leq\frac{3}{4},
\de
that is,
\be
\mP\left\{\sup_{0\leq s\leq t^*}|v(X_s,\mu_s)-v(x_0,\delta_{x_0})|\leq\epsilon_2\right\}\geq\frac{1}{4}.
\label{es1}
\ee

Now, set
\ce
S_1:&=&\inf\{t\geq0:v(X_t,\mu_t)<\varepsilon_1\},\\
S_2:&=&\inf\{t\geq S_1:v(X_t,\mu_t)>2\varepsilon_1\},\\
S_3:&=&\inf\{t\geq S_2:v(X_t,\mu_t)<\varepsilon_1\},\\
\vdots\\
S_{2k}:&=&\inf\{t\geq S_{2k-1}:v(X_t,\mu_t)>2\varepsilon_1\},\\
S_{2k+1}:&=&\inf\{t\geq S_{2k}:v(X_t,\mu_t)<\varepsilon_1\}, \quad k\in\mN^+,
\de
and then $\{S_i\}_{i=1}^\infty$ is a sequence of stopping times. Thus, by $(\ref{4.3})$ it holds that
\be
\infty>\mE\int_0^\infty v(X_s,\mu_s)ds&\geq&\sum_{k=1}^\infty\mE\left[I_{\{S_{2k}<\tau_n\}}\int_{S_{2k}}^{S_{2k+1}}v(X_s,\mu_s)ds\right]\no\\
&\geq&\varepsilon_1\sum_{k=1}^\infty\mE[I_{\{S_{2k}<\tau_n\}}(S_{2k+1}-S_{2k})]\no\\
&=&\varepsilon_1\sum_{k=1}^\infty\mE[I_{\{S_{2k}<\tau_n\}}\mE(S_{2k+1}-S_{2k}\mid\mathscr{G}_{S_{2k}})].\label{4.6}
\ee
We estimate $\mE(S_{2k+1}-S_{2k}\mid\mathscr{G}_{S_{2k}})$ on $\{S_{2k}<\tau_n\}$. First of all, by Lemma \ref{strmar} below, one can get that the strong solution $X_t$ of Eq.$(\ref{Eq2})$ have the strong Markov property for any $\{\mathscr{G}_t\}_{t\geq0}$-stopping time. Setting $\epsilon_2=\frac{\varepsilon_1}{2}$ and following the argument of Deng and Williams in \cite[P. 1241]{DW}, we obtain that on $\{S_{2k}<\tau_n\}$
\ce
\mE(S_{2k+1}-S_{2k}\mid\mathscr{G}_{S_{2k}})&\geq&\mE\((S_{2k+1}-S_{2k})I_{\{\sup\limits_{0\leq s\leq t^*}\mid v(\tilde{X}_s, \tilde{\mu}_s)-v(\tilde{X}_0, \tilde{\mu}_0)\mid\leq\frac{\varepsilon_1}{2}\}}\mid\mathscr{G}_{S_{2k}}\)\\
&\geq& t^*\mP\left\{{\sup_{0\leq s\leq t^*}|v(\tilde{X}_s, \tilde{\mu}_s)-v(\tilde{X}_0, \tilde{\mu}_0)|\leq\frac{\varepsilon_1}{2}}\mid\mathscr{G}_{S_{2k}}\right\}\\
&=& t^*\mP^{(S_{2k}, X_{S_{2k}})}\left\{\sup_{0\leq s\leq t^*}|v(X_s, \mu_s)-v(x_0, \delta_{x_0})|\leq\epsilon_2\right\}\\
&\overset{(\ref{es1})}\geq&\frac{t^*}{4},
\de
where $\tilde{X}_\cdot:=X_{\cdot+S_{2k}}$, $\tilde{\mu}_\cdot:=\mu_{\cdot+S_{2k}}$. Therefore, $(\ref{4.6})$ gives that
\ce
\frac{t^*}{4}\varepsilon_1\sum_{k=1}^\infty\mP\{S_{2k}<\tau_n\}<\infty.
\de
It follows from the Borel-Cantelli lemma that
\ce
\mP\{S_{2k}<\tau_n\ for\ infinitely\ many\ k\}=0.
\de
Note that
\ce
\{S_{2k}<\tau_n\ for\ infinitely\ many\ k\}&=&\{S_{2k}<\tau_n\ for\ infinitely\ many\ k\ and\ \tau_n=\infty\}\\
&&\cup\{S_{2k}<\tau_n\ for\ infinitely\ many\ k\ and\ \tau_n<\infty\}.
\de
Thus,
\ce
\mP\{S_{2k}<\tau_n\ for\ infinitely\ many\ k\ and\ \tau_n=\infty\}=0.
\de

Next, noting that $\sup\limits_{t\geq0}\mE(-v(X_t,\mu_t))^+=0$, by \cite[Theorem 3.15, P. 17]{IK} we get that $v(X_\infty, \mu_\infty)=\lim_{t\rightarrow\infty}v(X_t,\mu_t)$ exists and $\mE v(X_\infty,\mu_\infty)<\infty$. Therefore, it follows from the supermartingale inequality that
\ce
\mP\left\{\sup_{t\geq0}v(X_t,\mu_t)\geq\gamma_1(n)\right\}\leq\frac{\mE\({\sup\limits_{t\geq0}v(X_t,\mu_t)}\)}{\gamma_1(n)}\leq\frac{ v(x_0,\delta_{x_0})}{\gamma_1(n)},
\de
which gives
\ce
\mP\{\tau_n=\infty\}=\mP\left\{\sup_{t\geq0}| X_t|<n\right\}\geq\mP\left\{\sup_{t\geq0}v(X_t,\mu_t)<\gamma_1(n)\right\}\geq1-\frac{v(x_0,\delta_{x_0})}{\gamma_1(n)}.
\de
So, we have $\mP\{\tau_n=\infty\}\rightarrow1$ as $n\rightarrow\infty$, and furthermore
\ce
\mP\{S_{2k}<\infty\ for\ infinitely\ many\ k\}=0,
\de
which is a contradiction of $(\ref{4.4})$. This completes the proof.
\end{proof}

\bl\label{strmar}
Assume that $(\bf{H}_{1.1})$-$(\bf{H}_{1.2})$ hold, and $\tau$ is a $\{\mathscr{F}_t^W\}_{t\geq0}$-stopping time. Then the solution of Eq.$(\ref{Eq2})$ have the strong Markov property, that is
\be
\mE^{(0,x_0)}[h({\tau+t}, X_{\tau+t})|\mathscr{F}_\tau^W]=\mE^{(\tau, X_\tau)}[h(t, X_t)], \quad t\geq 0,
\label{strmarc}
\ee
where $h$ is a bounded Borel measurable function on $[0,\infty)\times\mR^d$.
\el
\begin{proof}
Let $X_t^{s,x_0}$ denote the unique strong solution of the following SDE:
$$
X_t^{s,x_0}=x_0+\int_s^t\tilde{b}(u, X^{s,x_0}_u)du+\int_s^t\tilde{\sigma}(u, X^{s,x_0}_u)d W_u, \quad t\geq s.
$$
So we have
\ce
X_{\tau+t}^{\tau,x_0}=x_0+\int_\tau^{\tau+t}\tilde{b}(u, X^{\tau,x_0}_u)du+\int_\tau^{\tau+t}\tilde{\sigma}(u, X^{\tau,x_0}_u)d W_u,
\de
and then
\be
X_{\tau+t}^{\tau,x_0}=x_0+\int_0^t\tilde{b}(\tau+r, X^{\tau,x_0}_{\tau+r})dr+\int_0^t\tilde{\sigma}(\tau+r, X^{\tau,x_0}_{\tau+r})d\tilde{W}_r,
\label{sto1}
\ee
where $\tilde{W}_{\cdot}:=W_{\tau+\cdot}-W_\tau$. By the strong Markov property of Brownian motion, it holds that $(\tilde{W}_s)_{s\geq0}$ is still a  Brownian motion, and $\{\mathscr{F}_s^{\tilde{W}}\}_{s\geq0}$ is independent of $\{\mathscr{F}_\tau^W\}$. Since $\{X_{\tau+t}^{\tau,x_0}\}$ is the unique strong solution of Eq.(\ref{sto1}), $\{X_{\tau+t}^{\tau,x_0}\}$ is adapted to $\{\mathscr{F}_s^{\tilde{W}}\}_{s\geq0}$ and then independent of $\{\mathscr{F}_\tau^W\}$. 

Besides, \cite[Proposition 3.7]{DH} shows that the pathwise uniqueness holds for Eq.$(\ref{sto1})$, which together with \cite[Proposition 3.20, P. 309]{IK} yields that the solution of Eq.$(\ref{sto1})$ is unique in law. So, $\{X_{\tau+t}^{\tau,x_0}\}_{t\geq 0}$ and $\{X_t^{0,x_0}\}_{t\geq 0}$ have the same distribution. 

Next, set 
$$
H(x_0, s, t,\omega):=X_t^{s, x_0}, \quad t\geq s\geq 0,
$$
and then by simple calculation, it holds that 
\ce
&&H(x_0,0,t,\omega)=X_t^{0,x_0}=X_t,\\
&&H(x_0,0,\tau+t,\omega)=H(X_\tau,\tau,\tau+t,\omega).
\de
Again set 
$$
\Psi(x_0, s, t, \omega):=h(t, H(x_0, s, t,\omega)),
$$
and then $\Psi$ is Borel measurable. Thus we can approximate $\Psi$ pointwise boundedly by functions of the form
$$
\sum_{k=1}^m \phi_k(x_0)\psi_k(s, t,\omega).
$$
So, we only prove $\Psi(x_0, s, t, \omega)=\sum_{k=1}^m \phi_k(x_0)\psi_k(s, t,\omega)$, and then take the limit to get (\ref{strmarc}) by the property of conditional expectations.

Finally, combining the above results, we have
\ce
&&\mE^{(0,x_0)}[h({\tau+t}, X_{\tau+t})|\mathscr{F}_\tau^W]=\mE^{(0,x_0)}[h({\tau+t}, X^{0,x_0}_{\tau+t})|\mathscr{F}_\tau^W]\\
&=&\mE^{(0,x_0)}[h({\tau+t}, H(x_0, 0, \tau+t,\omega))|\mathscr{F}_\tau^W]=\mE^{(0,x_0)}[h({\tau+t}, H(X_\tau, \tau, \tau+t,\omega))|\mathscr{F}_\tau^W]\\
&=&\mE^{(0,x_0)}[\Psi(X_\tau, \tau, \tau+t,\omega)|\mathscr{F}_\tau^W]=\mE^{(0,x_0)}\left[\sum_{k=1}^m \phi_k(X_\tau)\psi_k(\tau, \tau+t,\omega)|\mathscr{F}_\tau^W\right]\\
&=&\mE^{(0,x_0)}\left[\sum_{k=1}^m \phi_k(x_0)\psi_k(\tau, \tau+t,\omega)|\mathscr{F}_\tau^W\right]\bigg{|}_{x_0=X_\tau}=\mE^{(0,x_0)}\left[\Psi(x_0, \tau, \tau+t,\omega)|\mathscr{F}_\tau^W\right]\big{|}_{x_0=X_\tau}\\
&=&\mE^{(0,x_0)}\left[h({\tau+t}, H(x_0, \tau, \tau+t,\omega))|\mathscr{F}_\tau^W\right]\big{|}_{x_0=X_\tau}=\mE^{(0,x_0)}\left[h({\tau+t}, X^{\tau, x_0}_{\tau+t})|\mathscr{F}_\tau^W\right]\big{|}_{x_0=X_\tau}\\
&=&\mE^{(0,x_0)}\left[h({0+t}, X^{\tau, x_0}_{\tau+t})\right]\big{|}_{0=\tau,x_0=X_\tau}=\mE^{(0,x_0)}\left[h(t,X_t^{0,x_0})\right]\big{|}_{0=\tau,x_0=X_\tau}=\mE^{(\tau,X_\tau)}\left[h(t,X_t)\right].
\de
The proof is complete.
\end{proof}

\section{An example}\label{eg}

Now let us present an example to explain our results. 

\bx\label{e1}
Consider the following stochastic differential equation
\be\left\{\begin{array}{ll}
dX_t=-\(\int_{\mR}(X_t-my)\mu_t(dy)\)dt+\sum\limits_{k=1}^l\frac{\sin (kX_t)}{k^{3/2}}dW^k_t,\\
X_0=\xi,
\end{array}
\label{eqzui}
\right.
\ee
where $\xi$ is a $\sF_0$-measurable Gaussian random variable and $m>0$. It is easy to see that $b(x,\mu)=-x+m\int_{\mR}y\mu(\dif y)$ and $\sigma(x,\mu)=(1^{-\frac{3}{2}}\sin x,2^{-\frac{3}{2}}\sin 2 x,\ldots,l^{-\frac{3}{2}}\sin l x)$. So, one can justify that for $(x,\mu)\in\mR\times{\cM_{\lambda^2}(\mR)}$
\ce
|b(x,\mu)|^2&=&\(\int_{\mR}(x-my)\mu(dy)\)^2
\leq2\(|x|^2+m^2\left|\int_{\mR}y\mu(dy)\right|^2\)\\
&\leq&2\(|x|^2+m^2\int_{\mR}|y|^2\mu(dy)\)
\leq2(m^2+1)\(1+|x|^2+\int_{\mR}(1+|y|)^2\mu(dy)\)\\
&=&2(m^2+1)\(1+|x|^2+\|\mu\|_{\lambda^2}^2\),
\de
and
\ce
\|\sigma(x,\mu)\|^2&=&\sum_{k=1}^l\frac{|\sin kx|^2}{k^3}\
\leq\sum_{k=1}^l\frac{|\sin kx|}{k^3}\leq\sum_{k=1}^l\frac{1}{k^2}|x|
\leq S|x|\leq\frac{S}{2}(1+|x|^2),
\de
where $S:=\sum\limits_{k=1}^\infty\frac{1}{k^2}$. Thus, we have
\ce
|b(x,\mu)|^2+\|\sigma(x,\mu)\|^2\leq(2+\frac{S}{2}+2m^2)\(1+|x|^2+\|\mu\|_{\lambda^2}^2\).
\de

Moreover,  it holds that for $(x_1,\mu_1), (x_2,\mu_2)\in\mR^{d}\times{\cM_{\lambda^2}(\mR^d)}$, 
\ce
&&2\langle{x_1-x_2,b(x_1,\mu_1)-b(x_2,\mu_2)}\rangle+\|{\sigma(x_1,\mu_1)-\sigma(x_2,\mu_2)}\|^2\\
&=&2\langle{x_1-x_2,-\int_{\mR}(x_1-my)\mu_1(dy)+\int_{\mR}(x_2-my)\mu_2(dy)}\rangle\\
&&+\sum_{k=1}^l\frac{|\sin kx_1-\sin kx_2|^2}{k^3}\\
&\leq&2|x_1-x_2|\cdot\left|-\int_{\mR}(x_1-my)\mu_1(dy)+\int_{\mR}(x_2-my)\mu_2(dy)\right|\\
&&+4\sum_{k=1}^l\frac{|\sin \frac{k x_1- k x_2}{2}|^2}{k^3}\\
&\leq&|x_1-x_2|^2+\left|\int_{\mR}(x_1-my)\mu_1(dy)-\int_{\mR}(x_2-my)\mu_2(dy)\right|^2\\
&&+4\sum_{k=1}^\infty\frac{|\sin \frac{k x_1- k x_2}{2}|^2}{k^3}\\
&\leq&3|x_1-x_2|^2+2m^2\left|\int_{\mR}y\mu_1(dy)-\int_{\mR}y\mu_2(dy)\right|^2+C|x_1-x_2|^2\tilde{\kappa}(|x_1-x_2|)\\
&\leq&(3+2m^2+C)\(|x_1-x_2|^2+|x_1-x_2|^2\tilde{\kappa}(|x_1-x_2|)+\rho^2(\mu_1,\mu_2)\),
\de
where for $0<\eta<1/e$,
\ce
\tilde{\kappa}(x):=\left\{
\begin{array}{lcl}
\log x^{-1},&& 0<x\leq\eta,\\
\left(\left(\left(\log\eta^{-1}\right)^{\frac{1}{2}}-\frac{1}{2}\left(\log\eta^{-1}\right)^{-\frac{1}{2}}\right)x
+\frac{1}{2}\left(\log\eta^{-1}\right)^{-\frac{1}{2}}\eta\right)^2/x^2, &&x>\eta.
\end{array}
\right.
\de
Set $\kappa(x):=x^2+x^2\tilde{\kappa}(x)$, and then it is easy to justify that $\kappa(x)$ satisfies the conditions in $(\bf{H}_{1.2})$. Thus, the coefficients $b(x,\mu)$ and $\sigma(x,\mu)$ satisfy assumptions $(\bf{H}_{1.1})$-$(\bf{H}_{1.2})$, which yields that the equation has a strong solution.

If $v$ is taken as
\ce
v(x,\mu):=\(\int_{\mR}(x-my)\mu(dy)\)^2,
\de
then $v$ satisfies the following 
\ce
\int_{\mR}\(\mathscr{L^\mu}v(x,\mu)+\alpha v(x,\mu)\)\mu(dx)\leq \sum\limits_{k=1}^\infty\frac{1}{k^3}.
\de
Indeed, it is easily seen that
\ce
\partial_xv(x,\mu)&=&2\int_{\mR}(x-my)\mu(dy),\quad\partial_x^2v(x,\mu)=2,\\
\partial_{\mu}v(x,\mu)(z)&=&-2m\int_{\mR}(x-my)\mu(dy),\quad \partial_z\partial_{\mu}v(x,\mu)(z)=0.
\de
Hence it holds that
\ce
\mathscr{L^\mu}v(x,\mu)&=&-2\(\int_{\mR}(x-my)\mu(dy)\)^2+\sum_{k=1}^n\frac{|\sin kx|^2}{k^3}\no\\
&&+2m\int_{\mR}\[\int_{\mR}(z-my)\mu(dy)\cdot\int_{\mR}(x-my)\mu(dy)\]\mu(dz).
\de
Thus, integrating two sides of the above equality, we have
\ce
\int_{\mR}\(\mathscr{L^\mu}v(x,\mu)\)\mu(dx)&=&-2\int_{\mR}\(\int_{\mR}(x-my)\mu(dy)\)^2\mu(dx)+\int_{\mR}\sum_{k=1}^n\frac{|\sin kx|^2}{k^3}\mu(dx)+2m\cdot I,
\de
where $I:=\int_{\mR}\(\int_{\mR}\[\int_{\mR}(z-my)\mu(dy)\cdot\int_{\mR}(x-my)\mu(dy)\]\mu(dz)\)\mu(dx)$. Note that
\ce
I&=&\int_{\mR}\[\int_{\mR}\(\int_{\mR}(z-my)\mu(dy)\)\mu(dz)\cdot\int_{\mR}(x-my)\mu(dy)\]\mu(dx)\\
&=&\int_{\mR}\(\int_{\mR}(z-my)\mu(dy)\)\mu(dz)\cdot\int_{\mR}\(\int_{\mR}(x-my)\mu(dy)\)\mu(dx)\\
&=&\[\int_{\mR}\(\int_{\mR}(x-my)\mu(dy)\)\mu(dx)\]^2\\
&\leq&\int_{\mR}\(\int_{\mR}(x-my)\mu(dy)\)^2\mu(dx),
\de
where the last inequality is based on the H\"older inequality. Thus
\ce
\int_{\mR}\(\mathscr{L^\mu}v(x,\mu)\)\mu(dx)&\leq&-2\int_{\mR}\(\int_{\mR}(x-my)\mu(dy)\)^2\mu(dx)+\sum\limits_{k=1}^\infty\frac{1}{k^3}\\
&&+2m\int_{\mR}\(\int_{\mR}(x-my)\mu(dy)\)^2\mu(dx)\\
&=&(-2+2m)\int_{\mR}\(\int_{\mR}(x-my)\mu(dy)\)^2\mu(dx)+\sum\limits_{k=1}^\infty\frac{1}{k^3},
\de
i.e.
\ce
\int_{\mR}\(\mathscr{L^\mu}v(x,\mu)+(2-2m)v(x,\mu)\)\mu(dx)\leq\sum\limits_{k=1}^\infty\frac{1}{k^3}.
\de

Next, we justify that $v(x,\mu)$ satisfies (iii) of $(\bf{H}_{2.2})$. On one hand, it is obvious that
\ce
\int_{\mR}v(x,\mu)\mu(dx)&=&\int_{\mR}\(x-m\int_{\mR}y\mu(dy)\)^2\mu(dx)\\
&\leq&\int_{\mR}\[2|x|^2+2m^2\(\int_{\mR}y\mu(dy)\)^2\]\mu(dx)\\
&\leq&2\int_{\mR}|x|^2\mu(dx)+2m^2\int_{\mR}|y|^2\mu(dy)\\
&=&(2+2m^2)\int_{\mR}|x|^2\mu(dx).
\de
On the other hand, we have
\ce
\int_{\mR}v(x,\mu)\mu(dx)&=&\int_{\mR}\[|x|^2-2mx\int_{\mR}y\mu(dy)+m^2\(\int_{\mR}y\mu(dy)\)^2\]\mu(dx)\\
&\geq&\int_{\mR}\[|x|^2-2mx\int_{\mR}y\mu(dy)-m^2\(\int_{\mR}y\mu(dy)\)^2\]\mu(dx)\\
&=&\int_{\mR}|x|^2\mu(dx)-2m\int_{\mR}x\mu(dx)\int_{\mR}y\mu(dy)-m^2\(\int_{\mR}y\mu(dy)\)^2\\
&\geq&(1-2m-m^2)\int_{\mR}|x|^2\mu(dx).
\de
Finally, it holds that 
\ce
(1-2m-m^2)\int_{\mR}|x|^2\mu(dx)\leq\int_{\mR}v(x,\mu)\mu(dx)\leq (2+2m^2)\int_{\mR}|x|^2\mu(dx).
\de

In a word, if we choose appropriate $m$, $(\bf{H}_{2.2})$ is satisfied. For example, take $m=\frac{1}{4}$, then $\alpha=\frac{3}{2}$, $a_1=\frac{7}{16}, a_2=\frac{17}{8}$ and $M_1=\sum\limits_{k=1}^\infty\frac{1}{k^3}, M_2=M_3=0$. And then by Theorem \ref{sp} we know that the solution of Eq.(\ref{eqzui}) has the exponentially 2-ultimate boundedness. 
\ex

\bigskip

\textbf{Acknowledgements:}

Two authors would like to thank Professor Xicheng Zhang and Feng-Yu Wang for their valuable discussions. And they would also wish to thank two anonymous referees and Associate Editor for giving useful suggestions to improve this paper.

\end{document}